\documentclass[10pt,a4paper]{article}
\usepackage[utf8]{inputenc}
\usepackage[mathscr]{euscript}
\usepackage{amsfonts}
\usepackage{amsthm}
\usepackage{amsmath}
\usepackage{amssymb}
\usepackage{enumerate}

\title{A note on Fourier multipliers on the product Hardy spaces}
\author{Aleksander Pawlewicz}
\date{March 2022}

\theoremstyle{plain}
\newtheorem{theorem}{Theorem}
\newtheorem*{theorem*}{Theorem}
\newtheorem{lemma}[theorem]{Lemma}
\newtheorem{conjecture}[theorem]{Conjecture}

\theoremstyle{remark}

\begin{document}

\maketitle

\begin{abstract}
In this article, we present some results about Fourier multipliers on Hardy spaces in the product case. We mainly give descriptions of multipliers from the space $H^1(\mathbb{T}\times\mathbb{T})$ into the space $\ell^2$ and from the space $H^1(\mathbb{R}_+^2\times\mathbb{R}_+^2)$ into the space $L^1$.
\end{abstract}

\begin{section}*{Introduction}
In the article, we describe the basic facts of the theory of Fourier multipliers on product Hardy spaces. The two theorems we present provide characterization of such multipliers from the space  $H^1\left(\mathbb{T}\times\mathbb{T}\right)$ into $\ell^2$ and from the space $H^1\left(\mathbb{R}_+^2\times\mathbb{R}_+^2\right)$ into $L^1$. Because of the difficulties connected with the theory of product Hardy spaces and our imperfection, we decided to apply two different approaches to the two above-mentioned categories of multipliers. 

For the Fourier multipliers from $H^1\left(\mathbb{T}\times\mathbb{T}\right)$ into $\ell^2$ we base our considerations on the classical theory of Hardy spaces. To be more precise, we generalise ideas described in Duren's book \cite{Du} about ordinary one dimensional Hardy spaces defined on torus group $\mathbb{T}$. Definitions and results connected with this approach are presented in sections 1 and 2.

The Fourier multipliers from $H^1\left(\mathbb{R}_+^2\times\mathbb{R}_+^2\right)$ into $L^1$ are described in terms of the theory of product Hardy spaces developed mainly at the end of the 70's and in the 80's by mathematicians such as Sun-Young A. Chang, Robert Fefferman, Elias M. Stein and Richard F. Gundy. For more information about this type of Hardy spaces see \cite{ChFe3}, an excellent review on this subject from 1984. A multiplier results connected with this theory are presented in section 4.

\smallskip
In the case of one dimensional Hardy spaces
$$H^p(\mathbb{T})=\{f:\mathbb{T}\rightarrow\mathbb{C}:| f\in L^p(\mathbb{T}), \widehat{f}(n)=0 \mbox{ for } n<0\},$$
where $p\in[1,\infty)$, the problem of characterizing the Fourier multipliers operators dates back to the first half of the twentieth century. Hardy's inequality from 1927, see \cite[page 202]{HaLi}, can be considered the first result in this direction. It states that the sequence $(1/n)$ is a Fourier multiplier form $H^1(\mathbb{T})$ to $\ell^1$. 

The characterisation of Fourier multipliers from $H^1(\mathbb{T})$ into $\ell^2$ was given in paper \cite{DuSh}, while the description of Fourier multipliers from $H^1\left(\mathbb{R}^d\right)$ (in particular $H^1\left(\mathbb{R}\right)$) into $L^1$ was given in paper \cite{SlSt} (see also \cite{SzWo}). Article \cite{Os} contains the extensive review of results related to Fourier multipliers on the classical one-dimensional Hardy spaces.

The article is organised as follows. In Section 1, we present basic fact about the space $H^1\left(\mathbb{T}\times\mathbb{T}\right)$. In Section 2, we give a characterisation of multipliers from the space $H^1\left(\mathbb{T}\times\mathbb{T}\right)$ into $\ell^2$. Section 3 is devoted to basic facts connected with the $BMO\left(\mathbb{R}_+^2\times\mathbb{R}_+^2\right)$ space including different but equivalent definitions of this space. At the end, in Section 4, we present results which describes Fourier multipliers from $H^1\left(\mathbb{R}_+^2\times\mathbb{R}_+^2\right)$ space into the space $L^1$. We also present how our method works in the case of already known results of Fourier multipliers from $H^1\left(\mathbb{R}^d\right)$ into $L^1$. Moreover, we formulate the conjecture which describes how the condition for the Fourier multipliers in the product Hardy space may look. This conjecture simplifies the formula given in Theorem \ref{main_theorem_3}.
\end{section}

\begin{section}{Definition and basic facts about $H^1\left(\mathbb{T}\times\mathbb{T}\right)$}
The key tool in our consideration will be the Fourier transform. Let $\mathbb{T}$ denote the interval $[0,2\pi)$. We define the $n$-th Fourier coefficient, $n=(n_1,n_2, ...,n_d)\in\mathbb{Z}^d$, of locally integrable function $f:\mathbb{T}^d\rightarrow\mathbb{C}$ by the formula
$$\widehat{f}(n)=\frac{1}{(2\pi)^d}\int_{\mathbb{T}^d}f(x)e^{-in\cdot x}\,dx.$$

The two dimensional Hardy space is defined as follows:
$$H^1(\mathbb{T}\times\mathbb{T})=\{f:\mathbb{T}\times\mathbb{T}\rightarrow\mathbb{C}:| f\in L^1(\mathbb{T}\times\mathbb{T}), \widehat{f}(m,n)=0 \mbox{ for } m<0, n<0\},$$
with the norm 
$$\lVert f\rVert_{H^1}=\lVert f\rVert_{L^1}$$
for $f\in H^1(\mathbb{T}\times\mathbb{T})$.

By the Fourier multiplier operator we will mean the linear bounded operator $\Lambda:H^1(\mathbb{T}\times\mathbb{T})\rightarrow\ell^q$, $q\in[1,\infty)$, such that
$$\Lambda(f)=\left(\lambda(m,n)\widehat{f}(m,n)\right)_{m,n\in\mathbb{N}}.$$
Because functions from Hardy spaces have negative Fourier coefficients equal zero, it is enough to define the above sequences $(\lambda(m,n))$ only for non-negative integers.

\smallskip
In the proof of Theorem \ref{main_theorem_1} we will need one auxiliary theorem. We will call the two dimensional sequence $\left(\lambda(m,n)\right)$, $\lambda(m,n)\in\{0,1\}$ for all $m,n\in\mathbb{N}$, \textit{two-dimensional lacunary sequence} if in every dyadic rectangle $\left[2^k,2^{k+1}\right)\times\left[2^l,2^{l+1}\right)$, $k,l\in\mathbb{N}$, there is at most one non-zero element of the sequence $(\lambda(m,n))$. Then we have
\begin{theorem}[{{\cite[Thorem 1]{Av}}}]
\label{Aver}
There exists a positive constant $C>0$ such that for every two-dimensional lacunary sequence $\left(\lambda(m,n)\right)$ we have
\begin{equation*}
\left(\sum_{m,n=0}^\infty
\left|\lambda(m,n)\widehat{f}(m,n)\right|^2
\right)^{1/2}
\leq
C\lVert f\rVert_{L^1}
\end{equation*}
for $f\in H^1(\mathbb{T}\times\mathbb{T})$.
\end{theorem}

We also need to define the Fej\'er kernel. Let us denote the Fej\'er kernel by $F_n(x)$,
$$F_n(x)=\frac{1}{n+1}\left(\frac{\sin{\frac{n+1}{2}}}{\sin(x/2)}\right)^2=\sum_{j=-n}^n\left(1-\frac{|j|}{n+1}\right)e^{ijx}.$$
Then we will denote
$$F_{k,l}(x,y)=F_k(x)\cdot F_l(y)$$
for $k,l\in\mathbb{N}$. We also have for such $k, l$,
\begin{equation}
\label{Fejer_norm}
\lVert F_{k,l}\rVert_{L^1}=1.
\end{equation}
For more information about the definition of Fej\'er kernel and its properties we recommend \cite[Chapter I, section 2.5]{Ka}.
\end{section}

\begin{section}{Multipliers from $H^1\left(\mathbb{T}\times\mathbb{T}\right)$ into $\ell^2$}

\begin{theorem}
\label{main_theorem_1}
A square summable sequence $\lambda$ on $\mathbb{N}^2$ is a (Fourier) multiplier from the space $H^1\left(\mathbb{T}\times\mathbb{T}\right)$ into the space $\ell^2$ if and only if there exists a constant $C>0$ such that for all natural numbers $m$ and $n$ we have
\begin{equation}
\label{twierdzenie_Rudin}
\sum_{\alpha=2^m}^{2^{m+1}-1}\sum_{\beta=2^n}^{2^{n+1}-1}|\lambda(\alpha,\beta)|^2
\leq
C.
\end{equation}
\end{theorem}

Necessity of condition \eqref{twierdzenie_Rudin} is proved in the same way as Theorem 1 from \cite{Ru}.

\begin{proof}[Proof of Theorem \ref{main_theorem_1}]
Let us assume that $\lambda$ is a (Fourier) multiplier from the space $H^1\left(\mathbb{T}\times\mathbb{T}\right)$ into the space $\ell^2$ but
$$
\sup_{m,n\in\mathbb{N}}
\sum_{\alpha=2^m}^{2^{m+1}-1}\sum_{\beta=2^n}^{2^{n+1}-1}|\lambda(\alpha,\beta)|^2
=\infty.$$
Then there exists a sequence $\left(m_k,n_k\right)_{k\in\mathbb{N}}$ such that
$$\sum_{\alpha=2^{m_k}}^{2^{m_k+1}-1}\sum_{\beta=2^{n_k}}^{2^{n_k+1}-1}|\lambda(\alpha,\beta)|^2
>k.$$
Let $M$ denote the Fourier multiplier operator from the space $H^1\left(\mathbb{T}\times\mathbb{T}\right)$ into the space $\ell^2$ corresponding to the sequence $\lambda$. Then
$$\left\lVert M\left(e^{i2^{m_k+1}x}e^{i2^{n_k+1}y}F_{2^{m_k+1},2^{n_k+1}}(x,y)\right)\right\rVert_{\ell^2}^2
\geq
\frac{k}{4}.$$
Thus, because of \eqref{Fejer_norm}, we have
$$\lVert M\rVert\geq\frac{\sqrt{k}}{2}$$
for every $k\in\mathbb{N}$. This contradicts the boundedness of the norm of Fourier multiplier $M$ given by the sequence $\lambda$.

\smallskip
For the proof of sufficiency of condition \eqref{twierdzenie_Rudin} we will denote 
$$D_k=\{2^k,2^k+1,...,2^{k+1}-1\}.$$
Then we have, for $f\in H^1\left(\mathbb{T}\times\mathbb{T}\right)$, the following estimate
\begin{equation}
\label{szacowanie_M}
\begin{split}
\sum_{k=1}^\infty\sum_{l=1}^\infty
&\left|\lambda(k,l)\widehat{f}(k,l)\right|^2
\leq
\sum_{k=0}^\infty\sum_{l=0}^\infty
\sup_{\substack{\alpha\in D_k \\
\beta\in D_l}}
\left|\widehat{f}(\alpha,\beta)\right|^2
\sum_{\alpha\in D_k}\sum_{\beta\in D_l}
\left|\lambda(\alpha,\beta)\right|^2 \\
&\leq
\sup_{k,l\in\mathbb{N}}
\left(
\sum_{\alpha\in D_k}\sum_{\beta\in D_l}
\left|\lambda(\alpha,\beta)\right|^2
\right)
\sum_{k=0}^\infty\sum_{l=0}^\infty
\sup_{\substack{\alpha\in D_k \\
\beta\in D_l}}
\left|\widehat{f}(\alpha,\beta)\right|^2.
\end{split}
\end{equation}
The first term on the right hand side is bounded by assumption \eqref{twierdzenie_Rudin} and the second term is bounded by Theorem \ref{Aver}. Moreover, the operators $M:H^1\left(\mathbb{T}\times\mathbb{T}\right)\rightarrow\ell^2$ defined by
$$M(f)
=
\left(\lambda(k,l)\widehat{f}(k,l)\right)_{k,l\in\mathbb{N}}$$
have a closed graph. By \eqref{szacowanie_M} and the closed graph theorem, the operator $M$ is a bounded operator.
\end{proof}

\end{section}

\begin{section}{Product Hardy space $H^1(\mathbb{R}^2_+\times\mathbb{R}^2_+)$ and BMO space $BMO(\mathbb{R}^2_+\times\mathbb{R}^2_+)$}

In order to introduce the basic definition, we need an auxiliary function $\psi$. We will define it in the same way as it was done in paper \cite{ChFe1} (see the beginning of the proof on page 184). Let $\psi$ be a $C^1$ function supported on $[-1,1]$, even, with mean value $0$. We normalize $\psi$ so that
$$\int_0^\infty\left|\widehat{\psi}(s)\right|^2\frac{ds}{s}=1.$$Also, for $y_1,y_2>0$ and $x_1,x_2\in\mathbb{R}$ we define
\begin{equation}
\label{Psi}   
\Psi_y(x)=\frac{1}{y_1y_2}\psi\left(\frac{x_1}{y_1}\right)\psi\left(\frac{x_2}{y_2}\right),
\end{equation}
where $y=(y_1,y_2)$ and $x=(x_1,x_2)$. Then we have, for a locally integrable function $f$, the following formula
$$f(x)=\iint_{(t,y)\in\mathbb{R}_+^2\times\mathbb{R}_+^2}f(t,y)\Psi_y(x-t)\frac{dt\,dy}{y_1y_2},$$
and for dyadic rectangle $R$ we define
$$f_R(x)=\iint_{(t,y)\in R_+}f(t,y)\Psi_y(x-t)\frac{dt\,dy}{y_1y_2},$$
where
$$R_+=\left\{(t,y)\in\mathbb{R}_+^2\times\mathbb{R}_+^2| t\in R=I\times J, \frac{|I|}{2}<y_1\leq|I|, \frac{|J|}{2}<y_2\leq|J|\right\}.$$
For the function $f$ on $\mathbb{R}^2$ we define its double $S$-function
\begin{equation*}
S(f)(x)=\sqrt{
\iint_{\Gamma(x)}\left|f(t,y)\right|^2\frac{dt\,dy}{y_1^2y_2^2}},
\end{equation*}
where $\Gamma(x)=\Gamma(x_1)\times\Gamma(x_2)$ and
$$\Gamma(x_i)=\left\{(t_i,y_i)\in\mathbb{R}_+^2:|x_i-t_i|<y_i\right\},$$
for $i=1,2$.
All these facts are described on pages 182, 184, and 185 of \cite{ChFe1}. See also paragraphs 2.1 and 2.5 of \cite{We}.

Product Hardy space $H^1\left(\mathbb{R}_+^2\times\mathbb{R}_+^2\right)$ is a Banach space of integrable functions such that the norm 
\begin{equation*}
\lVert f\rVert_{H^1}=\lVert S(f)\rVert_{L^1}
\end{equation*}
is finite. There are a few equivalent definitions of the norm on this space, see part II, section 3 of \cite{ChFe3}, paragraph 1 of \cite{ChFe2}, section 2.1 of \cite{We} or \cite{GuSt}.

With all the above definitions we can formulate a basic theorem in this theory.
\begin{theorem}[{\cite[page 188]{ChFe1}}]
\label{BMO}
Suppose $\phi\in L^2(\mathbb{R}^2)$ satisfies
$$\int\phi(x_1,x_2)\,dx_1=\int\phi(x_1,x_2)\,dx_2=0$$
for all $(x_1,x_2)\in\mathbb{R}^2$. Then the following conditions on $\phi$ are equivalent:
\begin{enumerate}[a)]
\item $\phi$ is in the dual of $H^1(\mathbb{R}^2_+\times\mathbb{R}^2_+)$;
\item there exists a constant $C>0$ such that for every open set $\Omega\subset\mathbb{R}^2$ we have
$$\frac{1}{|\Omega|}\left\lVert\sum_{R\subset\Omega}\phi_R\right\rVert_2^2<C,$$
where $R$ are dyadic rectangles;
\item there exists a constant $C>0$ such that for every open set $\Omega\subset\mathbb{R}^2$ we have
$$\frac{1}{|\Omega|}\sum_{R\subset\Omega}S_R^2(\phi)<C,$$
where $R$ are dyadic rectangles and 
\begin{equation}
\label{SR2}
S_R^2(\phi)=\iint_{R_+}|\phi(t,y)|^2\frac{dt\,dy}{y_1y_2};
\end{equation}
\item there exists a constant $C>0$ such that for every open set $\Omega\subset\mathbb{R}^2$, there exists a function $\widetilde{\phi_\Omega}$ such that
$$\frac{1}{|\Omega|}\int_\Omega\left|\phi(t)-\widetilde{\phi_\Omega}(t)\right|^2\leq C.$$
\end{enumerate}
\end{theorem}

\bigskip

A locally integrable function $\lambda$ on $\mathbb{R}^2$ is a (Fourier) multiplier from the space $H^1(\mathbb{R}_+^2\times\mathbb{R}_+^2)$ into the space $L^p\left(\mathbb{R}^2\right)$, for $p\in[1,\infty)$, if and only if, the operator
$$f\longmapsto\widehat{f}(x)\lambda(x)$$
is a bounded linear operator from the space $H^1(\mathbb{R}_+^2\times\mathbb{R}_+^2)$ into the space $L^p\left(\mathbb{R}^2\right)$.

At the end, as a final remark, we would like to indicate that the theory of product Hardy and $BMO$ spaces in the sens of the Chang-Fefferman can also be developed for functions defined on the bidisc. This approach has been used in article \cite{PiWa} which describes how to obtain $BMO$ function by averaging dyadic $BMO$ functions on the bidisc.
\end{section}

\begin{section}{Multipliers from $H^1(\mathbb{R}^2_+\times\mathbb{R}^2_+)$ into $L^1$}

Before presenting the results in the product case we will formulate a new necessary condition for function to be a Fourier multiplier from ordinary Hardy spaces $H^1(\mathbb{R}^d)$ into $L^1$. Proof of this condition presents the method which we use in product case.
Let, as in \cite{SlSt}, $Q_\alpha^\varepsilon$ denote the cube:
$$Q_a^\varepsilon=\left\{x\in\mathbb{R}^d: \varepsilon\alpha_j-\frac{\varepsilon}{2}\leq x_j<\varepsilon\alpha_j+\frac{\varepsilon}{2}\right\},$$
where $\alpha=(\alpha_1, \alpha_2, ...,\alpha_n)\in\mathbb{Z}^d$ and $\varepsilon>0$.
Also we will write
$$\mathbb{T}_\varepsilon^d=\left[-\frac{\pi}{\varepsilon},\frac{\pi}{\varepsilon}\right)^d$$
for $\varepsilon>0$.

In this part we will use the following definitions of Fourier transform. For integrable function $f:\mathbb{R}^d\rightarrow\mathbb{C}$ we define
\begin{equation}
\label{transform_F}
\mathcal{F}(f)(\xi)=\frac{1}{(2\pi)^d}\int_{\mathbb{R}^d} f(x) e^{-i x\cdot\xi}\, dx.
\end{equation}
The inverse Fourier transform we will denote by $\mathcal{F}^{-1}$.
Moreover, for locally integrable function $f:\mathbb{T}_\varepsilon^d\rightarrow\mathbb{C}$, $\varepsilon>0$, we denote
\begin{equation}
\label{transform}
\widehat{f}(\varepsilon\alpha)=\frac{1}{\left(\frac{2\pi}{\varepsilon}\right)^d}\int_{\mathbb{T}_\varepsilon^d} f(x) e^{-i\varepsilon x\cdot\alpha}\, dx,
\end{equation}
where $\alpha=(\alpha_1, \alpha_2, ...,\alpha_n)\in\mathbb{Z}^d$.
We would like to draw the reader attention to the fact that
$$\widehat{f}(\varepsilon\alpha)=\widehat{f\left(\frac{\cdot}{\varepsilon}\right)}(\alpha).$$
In this situation we also have the formula for inverse Fourier transform
$$f(x)=\sum_{\alpha\in\mathbb{Z}^d} \widehat{f}(\varepsilon\alpha)e^{i\varepsilon \alpha\cdot x}$$
for $f:\mathbb{T}_\varepsilon^d\rightarrow\mathbb{C}$ and $x\in\mathbb{T}_\varepsilon^d$.

The proof which we want to present actually shows that expression which could characterise the norm of Fourier multiplier is equivalent to the $BMO$ norm of some function closely related to a given Fourier multiplier. 
In order to do that we need some definitions.

A complex valued locally integrable function $f$ on $\mathbb{R}^d$ is of a bounded mean oscillation, $f\in BMO$, if $\lVert f\rVert_{BMO}<\infty$ where
$$\lVert f\rVert_{BMO}=\sup_Q\frac{1}{|Q|}\int_Q\left|f(x)-\frac{1}{|Q|}\int_Qf(y)\,dy\right|\,dx,$$
where the supremum is taken over all cubes $Q$ in $\mathbb{R}^d$.
By the classical John-Nirenberg inequality we have the following corollary (see \cite[Corollary 7.1.9]{Gr}).
\begin{theorem}
\label{JohnNirenberg}
For $1<p<\infty$ we have
$$\lVert f\rVert_{BMO}
\approx
\sup_Q\left(\frac{1}{|Q|}\int_Q\left|f(x)-\frac{1}{|Q|}\int_Qf(y)\,dy\right|^p\,dx\right)^{1/p},$$
where the supremum is taken over all cubes $Q$ in $\mathbb{R}^d$.
\end{theorem}

Now we can formulate the preliminary result. Theorem 1 from article \cite{SlSt} states
\begin{theorem*}
Let $\mu$ be a positive Borel measure on $\mathbb{R}^d\setminus\{0\}$. Then
$$\sup\int_{\mathbb{R}^d}|\mathcal{F}(f)|d\,\mu<\infty,$$
where the supremum is taken over all $f$ in $H^1\left(\mathbb{R}^d\right)$ of norm $1$ if and only if
\begin{equation}
\label{warunek}
\sup_{\varepsilon>0}\left(\sum_{\alpha\in\mathbb{Z}^d\setminus\{0\}}\mu(Q_\alpha^\varepsilon)^2\right)^{1/2}<\infty.
\end{equation}
Moreover, the corresponding suprema are equivalent.
\end{theorem*}

Sufficiency of condition \eqref{warunek} in \cite{SlSt} was proved by the atomic decomposition of Hardy space. We prove using different approach the following statement.
\begin{theorem}
\label{klasyczne}
Let $\lambda$ be a locally integrable function on $\mathbb{R}^d$. Then
\begin{equation}
\label{norm1}
\sup_{\varepsilon>0}\left(\sum_{\alpha\in\mathbb{Z}^d\setminus\{0\}}\left|\widehat{\lambda}(\varepsilon\alpha)\right|^2\right)^{1/2}
\lesssim
\sup\int_{\mathbb{R}^d}|\mathcal{F}(f)(x)\mathcal{F}(\lambda)(x)|\,dx,
\end{equation}
where the supremum on the right hand side is taken over all $f$ in $H^1\left(\mathbb{R}^d\right)$ of norm $1$.
\end{theorem}

\begin{proof}
Let us denote by $\mathcal{H}$ the set of functions from the space $H^1\left(\mathbb{R}^d\right)$ of norm $1$. Then we have by Plancherel identity
\begin{equation*}
\begin{split}
\int_{\mathbb{R}^d}|
\mathcal{F}(f)(x)\mathcal{F}(\lambda)(x)
|\,dx
&\geq
\left|
\int_{\mathbb{R}^d}\mathcal{F}(f)(x)\overline{\mathcal{F}(\lambda)(x)}\,dx
\right| \\
&=
\left|
\int_{\mathbb{R}^d}f(x)\overline{\lambda(x)}\,dx
\right| 
\end{split}
\end{equation*}
Taking the supremum of the above inequality over functions f from $H^1\left(\mathbb{R}^d\right)$ of norm $1$ we obtain
\begin{equation}
\label{rachunek0}
\sup_{f\in\mathcal{H}}\int_{\mathbb{R}^d}|\mathcal{F}(f)(x)\mathcal{F}(\lambda)(x)|\,dx
\geq
\lVert\overline{\lambda}\rVert_{BMO}
=
\lVert\lambda\rVert_{BMO}.
\end{equation}

Now we calculate the $BMO$ norm of $\lambda$. Let $$P=\mathbb{T}_\varepsilon^d+p$$
be a cube of $\mathbb{R}^d$ with center $p\in\mathbb{R}^d$ and with side length $\varepsilon>0$. Then by Plancherel identity we have
\begin{equation}
\label{rachunek1}
\begin{split}
\frac{1}{|P|}\int_P&\left|\lambda(x)
-
\frac{1}{|P|}\int_P
\lambda(y)\,dy\right|^2dx \\
&=
\frac{1}{\left(\frac{2\pi}{\varepsilon}\right)^d}\int_{\mathbb{T}_\varepsilon^d}\left|\lambda(x+p)
-
\frac{1}{\left(\frac{2\pi}{\varepsilon}\right)^d}\int_{\mathbb{T}_\varepsilon^d}\lambda(y+p)\,dy
\right|^2dx \\
&=
\frac{1}{(2\pi)^d}\int_{\mathbb{T}^d}\left|\lambda\left(\frac{x+p}{\varepsilon}\right)
-
\frac{1}{(2\pi)^d}\int_{\mathbb{T}^d}\lambda\left(\frac{y+p}{\varepsilon}\right)\,dy
\right|^2dx \\
&=
\frac{1}{(2\pi)^d}\int_{\mathbb{T}^d}\left|\sum_{\alpha\in\mathbb{Z}^d\setminus\{0\}}\left\{\lambda\left(\frac{\cdot+p}{\varepsilon}\right)\right\}\widehat{}\,(\alpha)e^{i\alpha x}\right|^2dx \\
&=
\sum_{\alpha\in\mathbb{Z}^d\setminus\{0\}}\left|\left\{\lambda\left(\frac{\cdot+p}{\varepsilon}\right)\right\}\widehat{}\,(\alpha)\right|^2 \\
&=
\sum_{\alpha\in\mathbb{Z}^d\setminus\{0\}}\left|\left\{\lambda\left(\frac{\cdot}{\varepsilon}\right)\right\}\widehat{}\,(\alpha)e^{i\alpha p}\right|^2 \\
&=
\sum_{\alpha\in\mathbb{Z}^d\setminus\{0\}}\left|\widehat{\lambda\left(\frac{\cdot}{\varepsilon}\right)}(\alpha)\right|^2 \\
&=
\sum_{\alpha\in\mathbb{Z}^d\setminus\{0\}}\left|\widehat{\lambda}(\varepsilon\alpha)\right|^2.
\end{split}
\end{equation}
Taking square root and then supremum over $\varepsilon>0$ and $p\in\mathbb{R}^d$ of \eqref{rachunek1} we obtain by Theorem \ref{JohnNirenberg}
\begin{equation}
\label{rachunek2}
\lVert\lambda\rVert_{BMO}
\gtrsim
\sup_{\varepsilon>0}
\sqrt{\sum_{\alpha\in\mathbb{Z}^d\setminus\{0\}}|\widehat{\lambda}(\varepsilon\alpha)|^2}.
\end{equation}
Combining \eqref{rachunek0} and \eqref{rachunek2} we obtain
\begin{equation*}
\begin{split}
\sup_{f\in\mathcal{H}}\int_{\mathbb{R}^d}|\mathcal{F}(f)(x)\mathcal{F}(\lambda)(x)|\,dx
\geq
\lVert\lambda\rVert_{BMO} 
\gtrsim
\sup_{\varepsilon>0}
\sqrt{\sum_{\alpha\in\mathbb{Z}^d\setminus\{0\}}|\widehat{\lambda}(\varepsilon\alpha)|^2}.
\end{split}
\end{equation*}
\end{proof}

In order to formulate the main theorem of this section we need to split the real plane $\mathbb{R}^2$ into the set of identical rectangles. Let us define
$$
Q_{k,l}^{A,B}
=
I_k^A   \times   J_l^B
=
\left[kA-\frac{A}{2},kA+\frac{A}{2}\right)
\times
\left[lB-\frac{B}{2}, lB+\frac{B}{2}\right),$$
for $k,l\in\mathbb{Z}$ and $A,B>0$. 
Notice that the area of $Q_{k,l}^{A,B}$ is equal to $AB$,  $$\left|Q_{k,l}^{A,B}\right|=AB.$$

For the proof of the next theorem we need to change a little definition of Fourier transform given by formula \eqref{transform}. Let $\varepsilon=(\varepsilon_1,\varepsilon_2)$, $\varepsilon_1, \varepsilon_2>0$ and $$\mathbb{T}_\varepsilon^2=\left[-\frac{\pi}{\varepsilon_1},\frac{\pi}{\varepsilon_1}\right)
\times
\left[-\frac{\pi}{\varepsilon_2},\frac{\pi}{\varepsilon_2}\right).
$$
Then for locally integrable function $f:\mathbb{T}_\varepsilon^2\rightarrow\mathbb{C}$ we denote
\begin{equation}
\label{transform2}
\widehat{f}(\varepsilon\alpha)
=
\widehat{f}\left(\varepsilon_1\alpha_1,\varepsilon_2\alpha_2\right)
=
\frac{1}{\frac{(2\pi)^2}{\varepsilon_1\varepsilon_2}}\int_{\mathbb{T}_\varepsilon^2} f(x) e^{-i x\cdot(\varepsilon_1\alpha_1,\varepsilon_2\alpha_2)}\, dx,
\end{equation}
where $\alpha=(\alpha_1, \alpha_2)\in\mathbb{Z}^2$ and $\varepsilon=(\varepsilon_1,\varepsilon_2)$.
In this situation we also have the formula for inverse Fourier transform
$$f(x)=\sum_{\alpha\in\mathbb{Z}^d} \widehat{f}(\varepsilon\alpha)e^{ix\cdot(\varepsilon_1\alpha_1,\varepsilon_2\alpha_2)}$$
for $f:\mathbb{T}_\varepsilon^2\rightarrow\mathbb{C}$ and $x\in\mathbb{T}_\varepsilon^2$.

Now we can formulate one of the main theorems of this article.

\begin{theorem}
\label{main_theorem_2}
Let $\lambda$ be a locally integrable function on $\mathbb{R}^d$. Then for every bounded open set $\Omega$ in $\mathbb{R}^2$ we have
\begin{equation*}
\begin{split}
\sup
&\int_{\mathbb{R}^d}|\mathcal{F}(f)(x)\mathcal{F}(\lambda)(x)|\,dx \\
&\geq
\lVert\lambda\rVert_{BMO} \\
&\gtrsim
\sqrt{\frac{1}{|\Omega|}\sum_{R\subset\Omega}|R|
\sum_{k,l\in\mathbb{Z}}
\int_{|I|/2}^{|I|}\int_{|J|/2}^{|J|}
\left|
\widehat{\lambda*\Psi_y}(Ak,Bl)
\right|^2
\,\frac{dy}{y_1y_2}},
\end{split}
\end{equation*}
where the supremum on the left hand side is taken over all $f$ in $H^1(\mathbb{R}^2_+\times\mathbb{R}^2_+)$ of norm $1$
and the sum on the right hand side is taken over all dyadic rectangles $R=I\times J$ contained in $\Omega$; $A=\frac{2\pi}{|I|}$, $B=\frac{2\pi}{|J|}$.
\end{theorem}

\begin{proof}
In the proof we will use the definition of the $BMO\left(\mathbb{R}_+^2\times\mathbb{R}_+^2\right)$ space given by condition c) of Theorem \ref{BMO}, that is a function $\phi$ is in $BMO\left(\mathbb{R}_+^2\times\mathbb{R}_+^2\right)$ if there exists a constant $C>0$ such that for every open set $\Omega\subset\mathbb{R}^2$ we have
$$\sqrt{\frac{1}{|\Omega|}\sum_{R\subset\Omega}S_R^2(\phi)}<C,$$
where $R$ are dyadic rectangles.

Let us denote, similar to the proof of Theorem \ref{klasyczne}, by $\mathcal{H}$ the set of functions from the space $H^1(\mathbb{R}^2_+\times\mathbb{R}^2_+)$ of norm at most $1$. Then we have by Plancherel identity
\begin{equation*}
\begin{split}
\int_{\mathbb{R}^d}|
\mathcal{F}(f)(x)\mathcal{F}(\lambda)(x)
|\,dx
&\geq
\left|
\int_{\mathbb{R}^d}\mathcal{F}(f)(x)\overline{\mathcal{F}(\lambda)(x)}\,dx
\right| \\
&=
\left|
\int_{\mathbb{R}^d}f(x)\overline{\lambda(x)}\,dx
\right|
\end{split}
\end{equation*}
Taking the supremum of the above inequality over functions f in $\mathcal{H}$ we obtain
\begin{equation}
\label{rachunek10}
\sup_{f\in\mathcal{H}}\int_{\mathbb{R}^d}|\mathcal{F}(f)(x)\mathcal{F}(\lambda)(x)|\,dx
\geq
\lVert\overline{\lambda}\rVert_{BMO}
=
\lVert\lambda\rVert_{BMO}.
\end{equation}

Now let us calculate the $BMO$ norm of $\lambda$. We denote by $R=I\times J$ the dyadic rectangle centered in $(r_1,r_2)$ such that
$$R=I\times J=\mathbb{T}_{(A,B)}^2+(r_1,r_2),$$
where $A|I|=2\pi$ and $B|J|=2\pi$. Then by Plancherel identity we have
\begin{equation}
\label{rachunek11}
\begin{split}
S_R^2(\lambda)
&=
\iint_{R_+}|\lambda(t,y)|^2\frac{dt\,dy}{y_1y_2} \\
&=
\int_{|I|/2}^{|I|}\int_{|J|/2}^{|J|}
\iint_R|\lambda(t,y)|^2\,dt\frac{dy}{y_1y_2} \\
&=
\frac{1}{AB}\int_{|I|/2}^{|I|}\int_{|J|/2}^{|J|}
\int_{-\pi}^\pi\int_{-\pi}^\pi
\left|
\lambda*\Psi_y\left(\frac{t_1-r_1}{A},\frac{t_2-r_2}{B}\right)
\right|^2\,dt\frac{dy}{y_1y_2} \\
&=
|R|\int_{|I|/2}^{|I|}\int_{|J|/2}^{|J|}
\sum_{k,l\in\mathbb{Z}^2}
\left|
\left\{\lambda*\Psi_y\left(\frac{\cdot-r_1}{A},\frac{\cdot-r_2}{B}\right)\right\}\widehat{}\,(k,l)
\right|^2\,dt\frac{dy}{y_1y_2} \\
&=
|R|\int_{|I|/2}^{|I|}\int_{|J|/2}^{|J|}
\sum_{k,l\in\mathbb{Z}^2}
\left|
\widehat{\lambda*\Psi_y}(Ak,Bl)
\right|^2
\,\frac{dy}{y_1y_2} \\
&=
|R|
\sum_{k,l\in\mathbb{Z}^2}
\int_{|I|/2}^{|I|}\int_{|J|/2}^{|J|}
\left|
\widehat{\lambda*\Psi_y}(Ak,Bl)
\right|^2
\,\frac{dy}{y_1y_2},
\end{split}
\end{equation}
because $\lambda(t,y)=\lambda*\Psi_y(t)$.

Combining \eqref{rachunek10} and \eqref{rachunek11} we obtain by the definition of the $BMO\left(\mathbb{R}_+^2\times\mathbb{R}_+^2\right)$ space the following estimate: for every bounded open set $\Omega$ we have
\begin{equation*}
\begin{split}
\sup_{f\in\mathcal{H}}
&\int_{\mathbb{R}^d}|\mathcal{F}(f)(x)\mathcal{F}(\lambda)(x)|\,dx \\
&\geq
\lVert\lambda\rVert_{BMO} \\
&\gtrsim
\sqrt{\frac{1}{|\Omega|}\sum_{R\subset\Omega}|R|
\sum_{k,l\in\mathbb{Z}}
\int_{|I|/2}^{|I|}\int_{|J|/2}^{|J|}
\left|
\widehat{\lambda*\Psi_y}(Ak,Bl)
\right|^2
\,\frac{dy}{y_1y_2}},
\end{split}
\end{equation*}
where the sum on the right is taken over all dyadic rectangles $R=I\times J$ contained in the open bounded set $\Omega$.
\end{proof}

In fact we can prove a little stronger statement then presented in Theorem \ref{main_theorem_2}.
\begin{theorem}
\label{main_theorem_3}
Let $\lambda$ be a locally integrable function on $\mathbb{R}^d$. Then we have
\begin{equation*}
\begin{split}
&\sup_f
\int_{\mathbb{R}^d}|\mathcal{F}(f)(x)\mathcal{F}(\lambda)(x)|\,dx \\
&\approx
\sup_g
\sup_{\Omega}
\sqrt{\frac{1}{|\Omega|}\sum_{R\subset\Omega}|R|
\sum_{k,l\in\mathbb{Z}^2}
\int_{|I|/2}^{|I|}\int_{|J|/2}^{|J|}
\left|\left\{
\lambda*\mathcal{F}^{-1}(g)*\Psi_y\right\}\,\widehat{}\,(Ak,Bl)
\right|^2
\,\frac{dy}{y_1y_2}},
\end{split}
\end{equation*}
where the symbol $\sup_f$ on the left hand side denotes the supremum over all $f$ in $H^1(\mathbb{R}^2_+\times\mathbb{R}^2_+)$ of norm not bigger than $1$. Symbol $\sup_g$ denotes supremum over all functions $g$ such that $\lVert g\rVert_{L^\infty}\leq 1$ and symbol $\sup_\Omega$ denotes supremum over all bounded open subsets of $\mathbb{R}^2$.
The sum on the right hand side is taken over all dyadic rectangles $R=I\times J$ contained in $\Omega$; $A=\frac{2\pi}{|I|}$, $B=\frac{2\pi}{|J|}$.
\end{theorem}
\begin{proof}
The proof with the minor changes is almost the same as the proof of Theorem \ref{main_theorem_2}. We only need to see that we have by Plancherel identity
\begin{equation*}
\begin{split}
\int_{\mathbb{R}^d}|
\mathcal{F}(f)(x)\mathcal{F}(\lambda)(x)
|\,dx
&=
\int_{\mathbb{R}^d}
\left|
\mathcal{F}(f)(x)\overline{\mathcal{F}(\lambda)(x)}
\right|\,dx \\
&=
\sup_{\lVert g\rVert_{L^\infty\leq 1}}
\int_{\mathbb{R}^d}
\mathcal{F}(f)(x)\overline{\mathcal{F}(\lambda)(x)g(x)}
\,dx \\
&=
\sup_{\lVert g\rVert_{L^\infty\leq 1}}
\int_{\mathbb{R}^d}
\mathcal{F}(f)(x)\overline{\mathcal{F}\left(\lambda*\mathcal{F}^{-1}g(x)\right)}
\,dx \\
&=
\sup_{\lVert g\rVert_{L^\infty\leq 1}}
\int_{\mathbb{R}^d}f(x)\overline{\lambda*\mathcal{F}^{-1}g(x)}\,dx.
\end{split}
\end{equation*}
The rest of the proof is the same as in the proof of Theorem \ref{main_theorem_2}. We need only put in every place the expression $\lambda*\mathcal{F}^{-1}g$ instead of $\lambda$.
\end{proof}

At the end of this section we would like to state some conjecture. Condition given in Theorem \ref{main_theorem_3} can be probably simplified. It is eager to write the value
$$\left\{
\lambda*\mathcal{F}^{-1}(g)*\Psi_y\right\}\,\widehat{}\,(Ak,Bl)$$
from this condition in the following way
$$\widehat{\lambda}(Ak,Bl)g(Ak,Bl)\mathcal{F}(\Psi_y)(Ak,Bl),$$
but the situation is not so simple because of the use of different Fourier transform definitions and because of the problems with the supports of functions involved in this expression (function $\lambda$ is not periodic).
Although, using the following lemma we hope that condition form Theorem \ref{main_theorem_3} might be written in a better way.

\begin{lemma}
\label{lemat_prosty}
Let the function $\Psi$ be as it is define by \eqref{Psi}. Then for dyadic intervals $I$ and $J$ we have
\begin{equation}
\label{lemat_prosty_nier}
\begin{split}
\int_{|I|/2}^{|I|}\int_{|J|/2}^{|J|}\left|\mathcal{F}(\Psi_y)(Ak,Bl)\right|^2&\frac{dy_2\,dy_1}{y_1y_2} \\
&=\int_\frac{2\pi k}{2}^{2\pi k}\int_\frac{2\pi l}{2}^{2\pi l}\left|\mathcal{F}(\Psi)(y_1,y_2)\right|^2\frac{dy_2\,dy_1}{y_1y_2},
\end{split}
\end{equation}
where $\mathcal{F}(\Psi)(y_1,y_2)=\mathcal{F}(\Psi_{(1,1)})(y_1,y_2)$ and $k,l\in\mathbb{Z}$.
\end{lemma}

\begin{proof}
We have 
\begin{equation*}
\begin{split}
\int_{|I|/2}^{|I|}\int_{|J|/2}^{|J|}&\left|\mathcal{F}(\Psi_y)(Ak,Bl)\right|^2\frac{dy_2\,dy_1}{y_1y_2} \\
&=
\int_{|I|/2}^{|I|}\int_{|J|/2}^{|J|}\left|
\frac{1}{y_1y_2}
\mathcal{F}\left(\psi\left(\frac{\cdot}{y_1}\right)\right)(Ak)
\mathcal{F}\left(\psi\left(\frac{\cdot}{y_2}\right)\right)(Bl)\right|^2\frac{dy_2\,dy_1}{y_1y_2},
\end{split}
\end{equation*}
and by the Fubini theorem it is only enough to consider the integral
\begin{equation}
\label{Fubini}
\int_{|I|/2}^{|I|}
\left|
\frac{1}{y_1}
\mathcal{F}\left(\psi\left(\frac{\cdot}{y_1}\right)\right)(Ak)
\right|^2\frac{dy_1}{y_1}.
\end{equation}
By the definition of Fourier transform \eqref{transform} we have
\begin{equation}
\label{pomocniczy}
\frac{1}{y_1}
\mathcal{F}\left(\psi\left(\frac{\cdot}{y_1}\right)\right)(Ak)
=
\mathcal{F}(\psi)(Aky).
\end{equation}

Combining \eqref{Fubini} and \eqref{pomocniczy} we have by the change of variables
\begin{equation*}
\begin{split}
\int_{|I|/2}^{|I|}
\left|
\frac{1}{y_1}
\mathcal{F}\left(\psi\left(\frac{\cdot}{y_1}\right)\right)(Ak)
\right|^2\frac{dy_1}{y_1}
&=
\int_{|I|/2}^{|I|}
\left|
\mathcal{F}(\psi)(Aky_1)
\right|^2\frac{dy_1}{y_1} \\
&=
\int_{Ak|I|/2}^{Ak|I|}
\left|
\mathcal{F}(\psi)(y_1)
\right|^2\frac{dy_1}{y_1} \\
&=
\int_{\pi k}^{2\pi k}
\left|
\mathcal{F}(\psi)(y_1)
\right|^2\frac{dy_1}{y_1}.
\end{split}
\end{equation*}
This gives \eqref{lemat_prosty_nier}.
\end{proof}

Now we can formulate our conjecture.
\begin{conjecture}
Let $\lambda$ be a locally integrable function on $\mathbb{R}^d$. Then $\mathcal{F}(\lambda)$ is a (Fourier) multiplier from the space $H^1\left(\mathbb{R}_+^2\times\mathbb{R}_+^2\right)$ into the space $L^1$ if and only if there exists a constant $C>0$ such that for every bounded open set $\Omega$ we have
\begin{equation*}
\begin{split}
&\sqrt{\frac{1}{|\Omega|}\sum_{R\subset\Omega}|R|
\sum_{k,l\in\mathbb{Z}}
\left|
\widehat{\lambda}(Ak,Bl)
\right|^2
\int_\frac{2\pi k}{2}^{2\pi k}\int_\frac{2\pi l}{2}^{2\pi l}\left|\mathcal{F}(\Psi)(y_1,y_2)\right|^2\frac{dy_2\,dy_1}{y_1y_2}} \\
&=
\sqrt{\frac{1}{|\Omega|}\sum_{R\subset\Omega}|R|
\sum_{k,l\in\mathbb{Z}\setminus\{0\}}
\left|
\widehat{\lambda}(Ak,Bl)
\right|^2
\int_\frac{2\pi k}{2}^{2\pi k}\int_\frac{2\pi l}{2}^{2\pi l}\left|\mathcal{F}(\Psi)(y_1,y_2)\right|^2\frac{dy_2\,dy_1}{y_1y_2}}
<C,
\end{split}
\end{equation*}
where $R=I\times J$ are dyadic sub-rectangles of $\Omega$ and $A=\frac{2\pi}{|I|}$, $B=\frac{2\pi}{|J|}$.
\end{conjecture}

\end{section}

\bibliographystyle{plain}
\bibliography{sample}
\end{document}